\documentclass[12pt]{amsart}
\usepackage{latexsym, amsmath, amssymb}
\usepackage{epsfig}
\usepackage{amscd}
\usepackage{pstricks,pst-node,cite}

\newtheorem{thm}{Theorem}[section]
\newtheorem{lem}[thm]{Lemma}

\newtheorem{cor}[thm]{Corollary}

\theoremstyle{remark}
\newcommand{\C}[1]{\mathcal #1}
\newcommand{\B}[1]{\mathbb #1}

\def \Aut {\mbox{\rm Aut}}

\def \Inn {\mbox{\rm Inn}}

\ifx\blackandwhite\undefined
  \psset{linecolor=blue}\newrgbcolor{darkred}{.9 0 0}\newrgbcolor{emgreen}{0 .9 0}
  \newrgbcolor{lightgreen}{.7  .99 .7}\newrgbcolor{lightred}{.99 .7 .7} \else
  \psset{linecolor=blue}\newgray{darkred}{.7}\newgray{emgreen}{0.3}\newgray{pup}{.5}\newgray{xxx}{.5}
\fi
 \psset{linewidth=1pt,dash=3pt 3pt,doublesep=.05,dotsize=1pt 5}
\SpecialCoor

\begin{document}

\author{Dongseok Kim}
\address{Department of Mathematics \\Yeungnam University \\Kyongsan, 712-749, Korea} \email{dongseok@yu.ac.kr}
\thanks{}

\author{Jin Hwan Kim}
\address{Mathematics Education\\Yeungnam University \\Kyongsan, 712-749, Korea}
\email{kimjh@ynucc.yeungnam.ac.kr}

\author{Jaeun Lee}
\address{Department of Mathematics \\Yeungnam University \\Kyongsan, 712-749, Korea}
\email{julee@yu.ac.kr}
\thanks{The third author was supported by
Com$^2$MaC-KOSEF(R11-1999-054)}

\author{Dianjun Wang}
\address{Department of Mathematical Sciences\\Tsinghua University\\ Beijing, China}
\email{djwang@math.tsinghua.edu.cn}


\title[Enumerations of Cayley graphs]{Enumerations of Cayley graphs}
\begin{abstract}
We characterize the equivalence and the weak equivalence of Cayley
graphs for a finite group $\C{A}$. Using these characterizations, we
find enumeration formulae of the equivalence classes and weak
equivalence classes of Cayley graphs. As an application, we find the
number of weak equivalence classes of circulant graphs.
\end{abstract}

\maketitle

\section{Introduction}

Let $\C{A}$ be a finite group with identity $e$ and let $\Omega$ be
a set of generators for $\C{A}$ with the properties that $\Omega =
\Omega^{-1}$ and $e\not\in \Omega$, where $\Omega^{-1}=\{x^{-1}\,|\,
x\in \Omega\}$. The \emph{Cayley graph} $C(\C{A}, \Omega)$ is a
simple graph whose vertex-set and edge-set are defined as follows:
$$V(C(\C{A}, \Omega))=\C{A}~\mathrm{and}~ E(C(\C{A}, \Omega))=\{\{g, h\}\,
|\, g^{-1}h \in \Omega\}.$$ Because of its rich connections with
broad range of areas, it has been in the center of the research in
graph theory~\cite{BMPRS, DSS, PR, Rosenhouse}. Spectral estimations
of Cayley graphs have been studied~\cite{Cioaba, FMT}. It plays a
key role in the study of the geometry of hyperbolic
groups~\cite{kapo}. Recently, Li has found wonderful results on
edge-transitive Cayley graphs~\cite{Li, LL}.

The Cayley graph $C(\C{A}, \Omega)$ admits a natural $\C{A}$-action,
$\cdot: \C{A} \times C(\C{A}, \Omega) \to C(\C{A}, \Omega)$ defined
by $g\cdot g'=gg'$ for all $g,g'\in\C{A}$. A graph $G$ with an
$\C{A}$-action is called an \emph{$\C{A}$-graph}. So, every Cayley
graph $C(\C{A}, \Omega)$ is an $\C{A}$-graph. A graph isomorphism
$f:G\to H$ between two $\C{A}$-graphs is \emph{weak equivalence} if
there exists a group automorphism $\alpha:\C{A}\to \C{A}$ such that
$f(g\cdot u) = \alpha(g)\cdot f(u)$ for all $g\in \C{A}$ and $u\in
V(G)$. When $\alpha$ is the identity automorphism, we say that $f$
is an \emph{equivalence}. If there is a weak equivalence between
$\C{A}$-graphs $G$ and $H$, we say $G$ and $H$ are \emph{weak
equivalent}. Similarly, if there is an equivalence between
$\C{A}$-graphs $G$ and $H$, we say $G$ and $H$ are
\emph{equivalent}. For standard terms and notations, we refer
to~\cite{GTP1}.

Enumerations of the equivalence classes and weak equivalence classes
of some graphs have been studied~\cite{FKKL, KL1}. The purpose of
this article is to enumerate the equivalence classes and weak
equivalence classes of Cayley graphs for a finite group $\C{A}$.

The outline of this paper is as follows. In section~\ref{char}, we
characterize the weak equivalence of Cayley graphs for a finite
group $\C{A}$. Using these characterizations, we find enumeration
formulae of the equivalence classes and weak equivalence classes of
Cayley graphs in section~\ref{formula}. As an application, we find
the number of weak equivalence classes of circulant graphs in
section~\ref{appl}.

\section{A characterization of Cayley graphs}
\label{char}

Our definition of an weak equivalence between two Cayley graphs can
be interpolated as a color-consistence and direction preserving
graph isomorphism~\cite[Section 1.2.4]{GTP1}.

\begin{thm}\label{weakequiv}
Let  $C(\C{A}, \Omega)$ and $C(\C{A}, \Omega')$ be two Cayley
graphs. Then the followings are equivalent.
\begin{enumerate}
\item[{\rm (1)}] $C(\C{A}, \Omega)$ and $C(\C{A}, \Omega')$
           are weakly equivalent,
\item[{\rm (2)}] There exists  an isomorphism $\alpha:\C{A}\to
\C{A}$ such that $\alpha(\Omega)=\Omega'$.
\end{enumerate}
\end{thm}

\begin{proof}
$(1)\Rightarrow(2)$: Let $f:C(\C{A}, \Omega)\to C(\C{A}, \Omega')$
be a weak equivalence. Then there exists a group automorphism
$\tau:\C{A}\to \C{A}$ such that $f(g)=\tau(g)f(e)$ for each $g\in
\C{A}$. Let $x\in \Omega$. Then $\{e, x\}$ is an edge in $C(\C{A},
\Omega)$. Since $\{f(e), f(x)\}$ is an edge in
 $C(\C{A}, \Omega')$, $f(e)^{-1}f(x)=f(e)^{-1}\tau(x)f(e)$ is
 an element of $\Omega'$. Hence the map $\alpha:\C{A}\to\C{A}$
 defined by $\alpha(g)=f(e)^{-1}\tau(g)f(e)$ is a group
 isomorphism such that $\alpha(\Omega)=\Omega'$.

$(2)\Rightarrow(1)$: Let  $\alpha:\C{A}\to \C{A}$ be a group
automorphism such that $\alpha(\Omega)=\Omega'$. We define
$f:C(\C{A}, \Omega)\to C(\C{A}, \Omega')$ by $f(g)=\alpha(g)$. If
$\{g,h\}$ is an edge in $C(\C{A}, \Omega)$, then $g^{-1}h\in \Omega$
and
$f(g)^{-1}f(h)=\alpha(g)^{-1}\alpha(h)=\alpha(g^{-1}h)\in\alpha(\Omega)=\Omega'$.
Hence $f$ is a graph isomorphism such that
$f(gg')=\alpha(gg')=\alpha(g)\alpha(g')=\alpha(g)f(g')$, $i.e.$, $f$
is a weak equivalence.
\end{proof}

By using a similar method in Theorem~\ref{weakequiv}, we can have
the following theorem.

\begin{thm}\label{equiv}
Let  $C(\C{A}, \Omega)$ and $C(\C{A}, \Omega')$ be
two Cayley graphs. Then the followings are equivalent.
\begin{enumerate}
\item[{\rm (1)}] $C(\C{A}, \Omega)$ and $C(\C{A}, \Omega')$
           are  equivalent,
\item[{\rm (2)}] $\Omega$ and $\Omega'$ are conjugate in $\C{A}$,
$i.e.$, there exists an element $g\in \C{A}$ such that $g^{-1}\Omega
g=\Omega'$.
\end{enumerate}
\end{thm}

\section{Enumeration formulae}\label{formula}

For a finite group $\C{A}$, let  $$G_m(\C{A})=\{\Omega\subset
\C{A}\,:\, \Omega^{-1}=\Omega, <\Omega>=\C{A}, |\Omega|=m,e\not\in
\Omega \}.$$ Notice that $G_m(\C{A})$ contains all equivalence
classes of Cayley graphs $C(\C{A},\Omega)$ of degree $m$.

Let $H$ be a group of group automorphisms of $\C{A}$. $H$ admits a
natural action on $G_m(\C{A})$ by $\alpha\cdot
\Omega=\alpha(\Omega)$. By Theorem~\ref{weakequiv},
$\C{E}^w(\C{A},m)$, the number of weak equivalence classes of Cayley
graphs $C(\C{A},\Omega)$ of degree $m$, is equal to the number of
orbits of the $\Aut(\C{A})$ action on $G_m(\C{A})$, where
$\Aut(\C{A})$ is the group of all group isomorphisms of $\C{A}$.
Similarly, one can see that the number $\C{E}(\C{A},m)$ of
equivalence classes of Cayley graphs $C(\C{A},\Omega)$ of degree $m$
is equal to the number of orbits of the $\Inn(\C{A})$ action on
$G_m(\C{A})$ by Theorem~\ref{equiv}, where $\Inn(\C{A})$ is the
group of all inner automorphisms of $\C{A}$.

For any subset $S$ of $\C{A}$, let $O_2(S)=\{g\in S\,:\,g^2=e,
g\not=e\}$. We observe that
$$G_m(\C{A})=\bigcup_{k=0}^{\lfloor\frac{m}{2}\rfloor}
\{\Omega\in G_m(\C{A})\,:\, |O_2(\Omega)|=m-2k\}=:
\bigcup_{k=0}^{\lfloor\frac{m}{2}\rfloor}G_{m,k}(\C{A}).$$ This
implies that
$$\C{E}^w(\C{A},m)=\sum_{k=0}^{\lfloor\frac{m}{2}\rfloor}|G_{m,k}(\C{A})/\Aut(\C{A})|
\quad\mbox{ and }\quad
\C{E}(\C{A},m)=\sum_{k=0}^{\lfloor\frac{m}{2}\rfloor}|G_{m,k}(\C{A})/\Inn(\C{A})|.$$

Now, we aim to find a computational formula for the numbers
$|G_{m,k}(\C{A})/\Aut(\C{A})|$ and $|G_{m,k}(\C{A})/\Inn(\C{A})|$.
For each $x\in\C{A}$, let $\overline{x}=\{x,x^{-1}\}$. Let us denote
$\bar{G}_{m,k}(\C{A})$, the set of all $(m-k)$-tuples
$(\overline{x_1}, \overline{x_2},\ldots, \overline{x_k},
\overline{y_1},\ldots,\overline{y_{m-2k}})$ of distinct terms such
that
  $(1)$
$x_i^2\not=e, (i=1,2,\ldots, k)$, $(2)$ $y_j^{-1}=y_j\not=e,
(j=1,\ldots, m-2k)$, and $(3)$ $\{x_1, x_2,\ldots,
x_k,y_1,\ldots,y_{m-2k}\}$ generates $\C{A}$.

Let $S_k$ be the symmetric group on $k$ letters and let $\C{K}_k$ be
the direct product $S_k\times S_{m-2k}$ of $S_k$ and $S_{m-2k}$.
Define an $\Aut(\C{A})\times \C{K}_k$ action on
$\bar{G}_{m,k}(\C{A})$ by
\begin{align*}&(\alpha,\sigma,\tau)\cdot(\overline{x_1},\ldots, \overline{x_k},
\overline{y_1},\ldots,\overline{y_{m-2k}})\\&=(\overline{\alpha(x_{\sigma^{-1}(1)})},
\ldots, \overline{\alpha(x_{\sigma^{-1}(k)})},
\overline{\alpha(y_{\tau^{-1}(1)})},\ldots,\overline{\alpha(y_{\tau^{-1}(m-2k)})}).\end{align*}
Then it is not hard to show that
\begin{align*}|\bar{G}_{m,k}(\C{A})/\Aut(\C{A})\times
\C{K}_k|&=|G_{m,k}(\C{A})/\Aut(\C{A})|, \\
|\bar{G}_{m,k}(\C{A})/\Inn(\C{A})\times
\C{K}_k|&=|G_{m,k}(\C{A})/\Inn(\C{A})|.\end{align*} For a fixed
element $(\alpha, \sigma,\tau)$ in $\Aut(\C{A})\times \C{K}_k$, let
$F_{(\alpha,\sigma,\tau)}(m,k,\C{A})$ be the set all elements  in
$\bar{G}_{m,k}(\C{A})$  such that
$(\alpha,\sigma,\tau)\cdot(\overline{x_1},\ldots, \overline{x_k},
\overline{y_1},\ldots,\overline{y_{m-2k}})=(\overline{x_1},\ldots,
\overline{x_k}, \overline{y_1},\ldots,\overline{y_{m-2k}})$. The
following comes from the Burnside lemma.

\begin{thm}\label{countcal} Let $\C{A}$ be a finite group and $m$ be a
positive integer. Then the number of weak equivalence classes of
Cayley graphs $C(\C{A},\Omega)$ of degree $m$, $\C{E}^w(\C{A},m)$,
is
$$\sum_{k=0}^{\lfloor\frac{m}{2}\rfloor}|G_{m,k}(\C{A})/\Aut(\C{A})|
  =\sum_{k=0}^{\lfloor\frac{m}{2}\rfloor}\frac{\sum_{(\alpha, \sigma,
  \tau)\in \Aut(\C{A})\times
  \C{K}_k}|F_{(\alpha,\sigma,\tau)}(m,k,\C{A})|}{|\Aut(\C{A})|\,k!(m-2k)!},
  $$
and the number of equivalence classes of Cayley graphs
$C(\C{A},\Omega)$ of degree $m$, $\C{E}(\C{A},m)$, is
$$\sum_{k=0}^{\lfloor\frac{m}{2}\rfloor}|G_{m,k}(\C{A})/\Inn(\C{A})|
  =\sum_{k=0}^{\lfloor\frac{m}{2}\rfloor}\frac{\sum_{(\alpha, \sigma,
  \tau)\in \Inn(\C{A})\times
  \C{K}_k}|F_{(\alpha,\sigma,\tau)}(m,k,\C{A})|}{|\Inn(\C{A})|\,k!(m-2k)!}.
  $$
\end{thm}

\bigskip

Now, we will compute $|F_{(\alpha,\sigma,\tau)}(m,k,\C{A})|$. For
each subgroup $\C{S}$ of $\C{A}$ such that $\alpha(\C{S})=\C{S}$,
let us denote $\tilde{G}_{m,k}(\C{S})$ the set of all $(m-k)$-tuples
$(\overline{x_1}, \overline{x_2},\ldots, \overline{x_k},
\overline{y_1},\ldots,\overline{y_{m-2k}})$ of distinct elements in
$\C{S}$ such that $(1)$ $x_i^2\not=e, (i=1,2,\ldots, k)$, and $(2)$
$y_j^{-1}=y_j\not=e, (j=1,\ldots, m-2k)$. Then
$\tilde{G}_{m,k}(\C{S})$ is also an $\Aut(\C{A})\times \C{K}_k$ set.
Let $\tilde{F}_{(\alpha,\sigma,\tau)}(m,k,\C{S})$ be the set of all
elements in $\tilde{G}_{m,k}(\C{S})$ such that
$(\alpha,\sigma,\tau)\cdot(\overline{x_1},\ldots, \overline{x_k},
\overline{y_1},\ldots,\overline{y_{m-2k}})=(\overline{x_1},\ldots,
\overline{x_k}, \overline{y_1},\ldots,\overline{y_{m-2k}})$. It
follows from the M\"{o}bius inversion that
$$|F_{(\alpha,\sigma,\tau)}(m,k,\C{A})|=\sum_{\C{S}\le\C{A}, \alpha(\C{S})=\C{S}} \mu(\C{S})\,
 |\tilde{F}_{(\alpha,\sigma,\tau)}(m,k,\C{S})|,$$
where $\mu$ is the M\"obius function which assigns an integer
$\mu(\C{S})$ to each subgroup $\C{S}$ of $\C{A}$ such that
$\alpha(\C{S})=\C{S}$ by the recursive formula
$$\sum_{\C{S}'\geq \C{S}} \mu
(\C{S}')=\delta_{\C{S},\C{A}}=\left\{
\begin{array}{l}
1\ \mbox{if}\ \C{S}=\C{A},\\[.5ex] 0\ \mbox{if}\
\C{S} <\C{A}.
\end{array}
\right.$$

Let $\C{S}$ be a subgroup of $\C{A}$. We consider
$\Aut(\C{A})\times S_k$ action on the set $\tilde{G}_{2k,k}(\C{S})$
and $\Aut(\C{S})\times S_{m-2k}$ action on the set
$\tilde{G}_{m-2k,0}(\C{S})$. Then
$$|\tilde{F}_{(\alpha,\sigma,\tau)}(m,k,\C{S})|=
|\tilde{F}_{(\alpha,\sigma)}(2k,k,\C{S})|\,|\tilde{F}_{(\alpha,\tau)}(m-2k,0,\C{S})|,$$
for any subgroup $\C{S}$ of $\C{A}$. To complete the computation, we
need computational formulae for these two numbers
$|\tilde{F}_{(\alpha,\sigma)}(2k,k,\C{S})|$ and
$|\tilde{F}_{(\alpha,\tau)}(m-2k,0,\C{S})|$ for a subgroup $\C{S}$
of $\C{A}$ such that $\alpha(\C{S})=\C{S}$. For
$\alpha\in\Aut(\C{A})$ and a positive integer $n$, let
$$\tilde{F}_{(\alpha,n)}(\C{S})=\{g\in \C{S}\,:\, g\not=g^{-1}, \alpha^n(g)=g,
\alpha^l(g)\not=g~\mathrm{and}~ \alpha^l(g)\not=g^{-1}\, (l<n) \},$$
$$\tilde{I}_{(\alpha,n)}(\C{S})=\{g\in \C{A}\,:\,g\not=g^{-1}, \alpha^n(g)=g^{-1},
\alpha^l(g)\not=g~\mathrm{and}~ \alpha^l(g)\not=g^{-1}\, (l<n)\},$$
and
$$\tilde{F}^o_{(\alpha,n)}(\C{S})=\{g\in \C{S}\,:\, g^{-1}=g\not=e, \alpha^n(g)=g~\mathrm{and}~
\alpha^l(g)\not=g (l<n)\}.$$ For a fixed element $\sigma\in S_n$,
let $j_k(\sigma)$ be the number of disjoint $k$ cycles in the
factorization of $\sigma$ into disjoint cycles, $i.e.$,
$\sigma=\sigma_{j_1(\sigma)}\cdots\sigma_{j_n(\sigma)}$, where
$\sigma_{j_k(\sigma)}$ is the product of $j_k(\sigma)$ disjoint $k$
cycles.

\begin{lem}\label{counttil} Let
$\alpha\in\Aut(\C{A})$, $\sigma\in S_k$ and $\tau\in S_{m-2k}$.
Then for any subgroup $\C{S}$ of $\C{A}$ such that
$\alpha(\C{S})=\C{S}$, we have
$$|\tilde{F}_{(\alpha,\sigma,\tau)}(m,k,\C{S})|
= \prod_{r=1}^k
\prod_{s=0}^{j_r(\sigma)-1}
\left(\frac{|\tilde{F}_{(\alpha,r)}(\C{S})|+|\tilde{I}_{(\alpha,r)}(\C{S})|}{2}-rs\right)
\prod_{l=1}^{m-2k} \prod_{t=0}^{j_l(\tau)-1}
\left(|\tilde{F}^o_{(\alpha,l)}(\C{S})|-lt\right). $$
In particular,  if $|\C{A}|$ is odd, then we have
$$|\tilde{F}_{(\alpha,\sigma,\tau)}(m,k,\C{S})|
= \left\{\begin{array}{cl}
\displaystyle \prod_{r=1}^k
\prod_{s=0}^{j_r(\sigma)-1}
\left(\frac{|\tilde{F}_{(\alpha,r)}(\C{S})|+|\tilde{I}_{(\alpha,r)}(\C{S})|}{2}-rs\right)
& \mbox{\rm if $m$ is even and $k=\frac{m}{2}$,}\\[3ex]
0 &\mbox{\rm otherwise.}\end{array}\right.$$
\end{lem}

\begin{proof}
Notice that $(\overline{x_1},\ldots, \overline{x_k})$ is an element
of $\tilde{F}_{(\alpha,\sigma)}(2k,k,\C{S})$ if and only if
$\overline{x_i}=\overline{\alpha(x_{\sigma^{-1}(i)})}$ for each
$i=1,2,\ldots,k$. It means that if the length of the orbit of $i$
under the $<\sigma>$ is $r$, then $x_i=\alpha^{r}(x_i)$ or
$x_i^{-1}=\alpha^{r}(x_i)$, $i.e.$, $x_i\in
\tilde{F}_{(\alpha,r)}(\C{S})\cup \tilde{I}_{(\alpha,r)}(\C{S})$.
Since the number of orbits of length $r$ is  $j_r(\sigma)$, the set
$\{x_i,x_i^{-1},  \alpha(x_{\sigma(i)}),
\alpha(x_{\sigma(i)}^{-1}),\ldots, \alpha(x_{\sigma^{r-1}(i)}),
\alpha(x_{\sigma^{r-1}(i)}^{-1}) \}$ contains $2r$ numbers of
elements in $\tilde{F}_{(\alpha,r)}(\C{S})\cup
\tilde{I}_{(\alpha,r)}(\C{S})$, and $\{\overline{x_i},
\overline{\alpha(x_{\sigma(i)})},\ldots,
\overline{\alpha(x_{\sigma^{r-1}(i)})}\}=\{\overline{x_i^{-1}},
\overline{\alpha(x_{\sigma(i)}^{-1})},\ldots,
\overline{\alpha(x_{\sigma^{r-1}(i)}^{-1})}\}$. Since
$F_{(\alpha,r)}(\C{S})$ and $I_{(\alpha,r)}(\C{S})$ are disjoint, we
have
$$|\tilde{F}_{(\alpha,\sigma)}(2k,k,\C{S})|=\prod_{r=1}^k
\prod_{s=0}^{j_r(\sigma)-1}
\left(\frac{|\tilde{F}_{(\alpha,r)}(\C{S})|+|\tilde{I}_{(\alpha,r)}(\C{S})|}{2}-rs\right).$$
Similarly, we can show that
$$|\tilde{F}_{(\alpha,\tau)}(m-2k,0,\C{S})|=\prod_{l=1}^{m-2k} \prod_{t=0}^{j_l(\tau)-1}
\left(|\tilde{F}^o_{(\alpha,r)}(\C{S})|-lt\right).$$
Notice that if $|\C{A}|$ is odd, then $|\tilde{F}^o_{\alpha^r}(\C{S})|=0$.
 It completes the proof.
\end{proof}

We observe that if $\sigma_1$ and $\sigma_2$ are conjugate in $S_k$,
then
$|\tilde{F}_{(\alpha,\sigma_1)}(2k,k,\C{S})|=|\tilde{F}_{(\alpha,\sigma_2)}$$(2k$,
$ k$,$\C{S})|$ for any automorphism $\alpha$ of $\C{A}$ and any
subgroup $\C{S}$ of $\C{A}$ with $\alpha(\C{S})=\C{S}$. Similarly,
we can see that if $\tau_1$ and $\tau_2$ are conjugate in
$S_{m-2k}$, then
$|\tilde{F}^o_{(\alpha,\tau_1)}(m-2k,0,\C{S})|=|\tilde{F}_{(\alpha,\tau_2)}(m-2k,0,\C{S})|$
for any automorphism $\alpha$ of $\C{A}$ and any subgroup $\C{S}$ of
$\C{A}$ with $\alpha(\C{S})=\C{S}$. Moreover, the number of elements
in $S_m$ which are conjugate to $\sigma$ is equal to
$$\dfrac{m!}{j_1(\sigma)!2^{j_2(\sigma)}j_2(\sigma)!\cdots
m^{j_m(\sigma)}j_m(\sigma)!}.$$ Now, by the fact that $j_1(\sigma)+
2 j_2(\sigma)+\cdots+ m j_m(\sigma)=m$, Theorem~\ref{countcal} can
be reformulated as follows.

\begin{thm}\label{countcal2} Let $\C{A}$ be a finite group and $m$ be a
positive integer. Then we have
$$\begin{array}{l}
|\Aut(\C{A})|\,\,\C{E}^w(\C{A},m) \\
 =\displaystyle \sum_{k=0}^{\lfloor\frac{m}{2}\rfloor} \sum_{\alpha\in \Aut(\C{A})}
\sum_{\C{S}\le\C{A}, \alpha(\C{S})=\C{S}} \mu(\C{S})
\displaystyle 
\left(
\sum_{j_1+2j_2+\cdots+kj_k=k}\frac{\displaystyle
\prod_{r=1}^k
\prod_{t=0}^{j_r-1}
\left(\frac{|\tilde{F}_{\alpha^r}(\C{S})|+
|\tilde{I}_{\alpha^r}(\C{S})|}{2}-rs\right)}{
j_1!2^{j_2}j_2!\cdots k^{j_k}j_k!}\right)\\[3ex]
\displaystyle \hfill \times \left(
\sum_{j_1+2j_2+\cdots+(m-2k)j_{m-2k}=m-2k}\frac{\displaystyle
\prod_{l=1}^{m-2k}
\prod_{t=0}^{j_l-1}
\left(|\tilde{F}^o_{\alpha^l}(\C{S})|-lt\right)}{j_1!2^{j_2}j_2!\cdots (m-2k)^{j_{m-2k}}j_{m-2k}!}\right),
 \end{array} $$
and
$$\begin{array}{l}
|\Inn(\C{A})|\,\,\C{E}(\C{A},m) \\
 =\displaystyle \sum_{k=0}^{\lfloor\frac{m}{2}\rfloor} \sum_{\alpha\in \Inn(\C{A})}
\sum_{\C{S}\le\C{A}, \alpha(\C{S})=\C{S}} \mu(\C{S})
\displaystyle 
\left(
\sum_{j_1+2j_2+\cdots+kj_k=k}\frac{\displaystyle
\prod_{r=1}^k
\prod_{t=0}^{j_r-1}
\left(\frac{|\tilde{F}_{\alpha^r}(\C{S})|+
|\tilde{I}_{\alpha^r}(\C{S})|}{2}-rs\right)}{
j_1!2^{j_2}j_2!\cdots k^{j_k}j_k!}\right)\\[3ex]
\displaystyle \hfill \times \left(
\sum_{j_1+2j_2+\cdots+(m-2k)j_{m-2k}=m-2k}\frac{\displaystyle
\prod_{l=1}^{m-2k}
\prod_{t=0}^{j_l-1}
\left(|\tilde{F}^o_{\alpha^l}(\C{S})|-lt\right)}{j_1!2^{j_2}j_2!\cdots (m-2k)^{j_{m-2k}}j_{m-2k}!}\right).
 \end{array} $$
\end{thm}

Now, we will compute $|\tilde{F}_{(\alpha,r)}(\C{S})|$,
$|\tilde{I}_{(\alpha,r)}(\C{S})|$, and
$|\tilde{F}^o_{(\alpha,r)}(\C{S})|$. For convenience, let
\begin{align*} &\tilde{\tilde{F}}_{(\alpha,r)}(\C{S})=\{g\in \C{S}\,:\,
\alpha^r(g)=g, g^{-1}\not=g\},\\
&\tilde{\tilde{I}}_{(\alpha,r)}(\C{S})=\{g\in \C{S}\,:\,
\alpha^r(g)=g^{-1}, g^{-1}\not=g\},\\
&\tilde{\tilde{F}}^o_{(\alpha,r)}(\C{S})=\{g\in \C{S}\,:\,
\alpha^r(g)=g, g^{-1}=g\not=e\}.\end{align*}

\begin{lem} Let $\C{A}$ be a finite group and let
$\alpha$ be an automorphism on $\C{A}$ of order $|<\alpha>|$. Then
for any positive integer $r$ and any subgroup $\C{S}$ of $\C{A}$
such that $\alpha(\C{S})=\C{S}$, we have
$$|\tilde{F}_{\alpha^r}(\C{S})|=\left\{
\begin{array}{ll}\displaystyle
\sum_{d|r}\mu\left(\frac{r}{d}\right)\,\left|\tilde{\tilde{F}}_{\alpha^d}(\C{S})\right|
 & \mbox{\rm if $r$ is a divisor of $|\!<\alpha>\!|$ and odd,} \\[4ex]
\displaystyle\sum_{d|r}\mu\left(\frac{r}{d}\right)\,\left|\tilde{\tilde{F}}_{\alpha^d}(\C{S})\right|
-|\tilde{I}_{\alpha^\frac{r}{2}}(\C{S})|
      & \mbox{\rm if $r$ is a divisor of $|\!<\alpha>\!|$ and even, }\\[4ex]
0 & \mbox{\rm if $r$ is not a divisor of $|\!<\alpha>\!|$,}
\end{array}\right.$$
$$|\tilde{F}^o_{\alpha^r}(\C{S})|=\left\{\begin{array}{cl}\displaystyle
\sum_{d|r}\mu\left(\frac{r}{d}\right)\,\left|\tilde{\tilde{F}}^o_{\alpha^d}(\C{S})\right|
 & \mbox{\rm if
$r$ is a divisor of $|\!<\alpha>\!|$,}\\
0 & \mbox{\rm otherwise,}
\end{array}\right.$$
 and
 $$
|\tilde{I}_{\alpha^r}(\C{S})|=\left\{\begin{array}{cl}\displaystyle
\sum_{d|r ~and~\frac{r}{d}
~is~odd}\mu\left(\frac{r}{d}\right)\left|\tilde{\tilde{I}}_{\alpha^d}(\C{S})\right|
& \mbox{\rm if
$2r$ is a divisor of $|\!<\alpha>\!|$,}\\[1ex]
0 & \mbox{\rm otherwise.} \end{array}\right.
$$ \label{countttil}
\end{lem}

\begin{proof}
Let $x\in\C{S}$ such that $\alpha^r(x)=x$ and $x\not=x^{-1}$. If
$\alpha^l(x)=x$ for some $l$, then $\alpha^d(x)=x$, where $d=(r,l)$
is the greatest common divisor of $r$ and $l$. It implies that if
$\alpha^l(x)=x$ then $x\in\tilde{\tilde{F}}_{\alpha^d}(\C{S})$ for
some divisor $d$ of $r$. Since $\alpha^{(r, |\!<\alpha>\!|)}(x)=x$,
$|\tilde{F}_{\alpha^r}(\C{S})|\not=0$ when $r=(r,
|\!\!<\alpha>\!\!|)$, $i.e.$, $r$ is a divisor of
$|\!\!<\alpha>\!\!|$.

If $\alpha^l(x)=x^{-1}$ for $l<r$, then $\alpha^{2l}(x)=x$ and hence
$\alpha^{(r,2l)}(x)=x$. It implies that if $d=(r,2l)< r$ then $x\in
\tilde{\tilde{F}}_{\alpha^d}(\C{S})$, and if $d=(r,2l)= r$, then
$2l=r$ and $x\in \tilde{\tilde{I}}_{\alpha^\frac{r}{2}}(\C{S})$.
Now, we can see that
$$\tilde{F}_{\alpha^r}(\C{S})=\tilde{\tilde{F}}_{\alpha^r}(\C{S})-\left(
\bigcup_{d|r~and~d\not=r}\tilde{\tilde{F}}_{\alpha^d}(\C{S})\cup\tilde{\tilde{I}}_{\alpha^\frac{r}{2}}(\C{S})\right).
$$

Notice that if $\alpha^s(x)=x$ and $s|t$ then $\alpha^t(x)=x$,
$i.e.$, $\tilde{\tilde{F}}_{\alpha^s}(\C{S})\subset
\tilde{\tilde{F}}_{\alpha^t}(\C{S})$ for each $s$ and $t$ with
$s|t$. Since
$$\tilde{\tilde{I}}_{\alpha^\frac{r}{2}}(\C{S})-\bigcup_{d|r~and~d\not=r}\tilde{\tilde{F}}_{\alpha^d}(\C{S})
=\tilde{I}_{\alpha^\frac{r}{2}}(\C{S}),
$$
we have
$$|\tilde{F}_{\alpha^r}(\C{S})|=\left\{
\begin{array}{ll}\displaystyle
\sum_{d|r}\mu\left(\frac{r}{d}\right)\,\left|\tilde{\tilde{F}}_{\alpha^d}(\C{S})\right|
 & \mbox{if
$r$ is a divisor of $|\!<\alpha>\!|$ and odd,} \\[4ex]
\displaystyle\sum_{d|r}\mu\left(\frac{r}{d}\right)\,\left|\tilde{\tilde{F}}_{\alpha^d}(\C{S})\right|
-|\tilde{I}_{\alpha^\frac{r}{2}}(\C{S})|
      & \mbox{if
$r$ is a divisor of $|\!<\alpha>\!|$ and even,}\\[4ex]
0 & \mbox{if $r$ is not a divisor of $|\!<\alpha>\!|$.}
\end{array}\right.$$
By a method similar to the computation of
$|\tilde{F}_{\alpha^r}(\C{S})|$, we have
$$|\tilde{F}^o_{\alpha^r}(\C{S})|=
\sum_{d|r}\mu\left(\frac{r}{d}\right)\,|\tilde{\tilde{F}}^o_{\alpha^d}(\C{S})|.$$

Let $x\in\C{S}$ such that $\alpha^r(x)=x^{-1}$ and $x\not=x^{-1}$.
If $\alpha^l(x)=x$ for some $l$, then $\alpha^{(2r,l)}(x)=x$. In
particular,  $\alpha^{(2r,|<\alpha>|)}(x)=x$. This implies that
$|\tilde{I}_{\alpha^r}(\C{S})|\not=0$ when $2r=(2r, |<\alpha>|)$,
$i.e.$, $2r$ is a divisor of $|\!<\alpha>\!|$. For $l\le r$, Put
$d=(2r,l)$. Then $\alpha^d(x)=x$. If $d$ is a divisor of $r$, then
$x=\alpha^{d\,\frac{r}{d}}(x)=\alpha^r(x)=x^{-1}$. Since
$x\not=x^{-1}$, $d$ can not be a divisor of $r$. Since $d|2r$, $d$
is even and $\dfrac{2r}{d}$ is odd. Hence,
$x^{-1}=\alpha^r(x)=\alpha^{\frac{d}{2}\,\frac{2r}{d}}(x)=\alpha^{\frac{d}{2}}(x)$,
$i.e.$, $\alpha^{d'}(x)=x^{-1}$ for some $d'|r$ and $\dfrac{r}{d'}$
is odd. If $\alpha^l(x)=x^{-1}$ for some $l\le r$, then
$\alpha^{(2r,2l)}(x)=x$. Put $(2r,2l)=2(r,l)=2d$. Then
$\alpha^{2d}(x)=x$ and
$x^{-1}=\alpha^{r}(x)=\alpha^{d\,\frac{r}{d}}(x)$. Since
$x\not=x^{-1}$, $\dfrac{r}{d}$ is odd and $\alpha^d(x)=x^{-1}$. Now,
we can see that
$$\tilde{I}_{\alpha^r}(\C{S})=\displaystyle\tilde{\tilde{I}}_{\alpha^r}(\C{S})-\bigcup_{d|r,~
d\not=r~and~\frac{r}{d}~is
~odd}\tilde{\tilde{I}}_{\alpha^d}(\C{S}),$$ and hence,
$$|\tilde{I}_{\alpha^r}(\C{S})|=
\sum_{d|r ~and~\frac{r}{d}
~is~odd}\mu\left(\frac{r}{d}\right)\,|\tilde{\tilde{I}}_{\alpha^d}(\C{S})|.$$
It completes the proof.
\end{proof}

Let $\C{A}$ be a finite abelian group. Then
$\Inn(\C{A})=\{id_{\C{A}}\}$. Let $\C{S}$ be a subgroup of $\C{A}$.
Then
$|\tilde{F}_{1}(\C{S})|+|\tilde{I}_{1}(\C{S})|=|\C{S}|-|O_2(\C{S})|-1$
and $|\tilde{F}^o_{1}(\C{S})=|O_2(\C{S})|$. Now, by
Theorem~\ref{countcal2} and Lemma~\ref{countttil}, we have the
following corollary.

\begin{cor} Let $\C{A}$ be a finite abelian group and $m$ be a
positive integer. Then the number of weak equivalence classes of
Cayley graphs $C(\C{A},\Omega)$ of degree $m$, $\C{E}(\C{A},m)$ is
$$
\sum_{k=0}^{\lfloor\frac{m}{2}\rfloor}\sum_{\C{S}\le\C{A}}
\mu(\C{S}) \left(
\begin{matrix}\frac{1}{2}(|\C{S}|-|O_2(\C{S})|-1)\\
k\end{matrix}\right) \left(\begin{matrix}|O_2(\C{S})|\\
m-2k\end{matrix}\right),$$ where $O_2(\C{S})=\{g\in\C{S}~:~g^2=e,
g\not=e\}.$ In particular, if $\C{A}$ is odd, then
$$
\C{E}(\C{A},m) =\left\{\begin{array}{cl}
\displaystyle  \sum_{\C{S}\le\C{A}} \mu(\C{S})
\left(\begin{matrix}\frac{1}{2}(|\C{S}|-1)\\ \frac{m}{2}\end{matrix}\right) & \mbox{\rm if $m$ is even,}\\[2ex]
0 & \mbox{\rm if $m$ is odd.}\end{array} \right.$$
\label{countcaleabel}
\end{cor}


\section{Applications to circulant graphs}\label{appl}

Let $\B{Z}_n$ be the additive cyclic group of order $n$. A connected
\emph{circulant graph} is a Cayley graph $C(\B{Z}_n,\Omega)$ for the
cyclic group $\B{Z}_n$ of order $n$. Let $p$ be a prime number. Two
circulant graphs for a cyclic group $\B{Z}_p$ are isomorphic if and
only if they are weakly equivalent~\cite{ET70}. So, the number
$\C{E}^w(\B{Z}_p,m)$ of weak equivalence classes and the number of
isomorphism classes of Cayley graphs $C(\B{Z}_p,\Omega)$ of degree
$m$ are equal.

We identify $\Aut(\B{Z}_n)$ with the set of all elements of
$\B{Z}_n$ which are relatively prime to $n$, that is, the set
$\{\alpha\in\B{Z}_n\,:\, (\alpha,n)=1\}$. Notice that
$\Aut(\B{Z}_n)$ has $\phi(n)$ elements, where $\phi$ is the Euler
function. Notice that the number of elements $g$ in $\B{Z}_n$ such
that $g=-g$ is one if $n$ is odd or two if $n$ is even. Moreover,
such elements are fixed by every automorphism $\alpha$ of $\B{Z}_n$,
$i.e.$, $\alpha(g)=g$. Since $\alpha^r$ is also an automorphism for
any automorphism $\alpha$ and any integer $r$, we have
 $$|\tilde{F}^o_{\alpha^r}(\B{Z}_n)|=\left\{\begin{array}{ll}1
& \mbox{\rm if $r=1$ and $n$ is even,}\\[1ex]
0 & \mbox{\rm otherwise.}\end{array}\right.$$

Now, we aim to compute $|\tilde{F}_{\alpha^r}(\B{Z}_n)|$ and
$|\tilde{I}_{\alpha^r}(\B{Z}_n)|$. By Lemma~\ref{countttil}, it is
sufficient to compute $|\tilde{\tilde{F}}_{\alpha^r}(\B{Z}_n)|$ and
$|\tilde{\tilde{I}}_{\alpha^r}(\B{Z}_n)|$. Let $\alpha$ be an
automorphism and let $r$ be an integer. Then $g\in
\tilde{\tilde{F}}_{\alpha^r}(\B{Z}_n)$ if and only if
$\alpha^r(g)=g$ and $2g\not=0$, $i.e.$, $(\alpha^r-1)g=0$ and
$2g\not=0$, and $g\in \tilde{\tilde{I}}_{\alpha^r}(\B{Z}_n)$ if and
only if $\alpha^r(g)=-g$ and $2g\not=0$, $i.e.$, $(\alpha^r+1)g=0$
and $2g\not=0$. Hence we have the following lemma.

\begin{lem} Let $\alpha\in\B{Z}_n$ such that $(\alpha,n)=1$ and for any natural
number $r$, we have
$$|\tilde{\tilde{F}}_{\alpha^r}(\B{Z}_n)|=\left\{\begin{array}{ll}
(\alpha^r-1,n)-1 & \mbox{\rm if $n$ is odd,}\\[1ex]
(\alpha^r-1,n)-2 & \mbox{\rm if $n$ is even,}\end{array}\right.
$$ $$
|\tilde{\tilde{I}}_{\alpha^r}(\B{Z}_n)|=\left\{\begin{array}{ll}(\alpha^r+1,n)-1
& \mbox{\rm if $n$ is odd,}\\[1ex]
(\alpha^r+1,n)-2 & \mbox{\rm if $n$ is even.}\end{array}\right.$$
\label{countttilzn}
\end{lem}

\begin{cor}\label{counttilzp}
Let $p$ be an odd prime and let $\alpha\in\B{Z}_p$ such that
$(\alpha,p)=1$. Then for any natural number $r$, we have
$$|\tilde{F}_{\alpha^r}(\B{Z}_p)|=\left\{\begin{array}{ll}
p-1 & \mbox{\rm if $r$ is odd, the order of $\alpha$ is $r$ and $r|(p-1)$,}\\
0   & \mbox{\rm otherwise,}\end{array}\right.
$$
and
$$|\tilde{I}_{\alpha^r}(\B{Z}_p)|=\left\{\begin{array}{ll}p-1
& \mbox{\rm if order of $\alpha$ is $2r$ and $r|(\dfrac{p-1}{2})$,} \\
0& \mbox{\rm otherwise.} \end{array}\right.$$
\end{cor}

From Corollary~\ref{counttilzp}, we can see that
$$|\tilde{F}_{\alpha^r}(\B{Z}_p)|+ |\tilde{I}_{\alpha^r}(\B{Z}_p)|
=\left\{\begin{array}{ll}
p-1 & \mbox{\rm if $r$ is odd, order of $\alpha$ is $r$ and $r|(p-1)$,}\\[.8ex]
p-1 & \mbox{\rm if order of $\alpha$ is $2r$ and $r|(\dfrac{p-1}{2})$,}\\[1ex]
0   & \mbox{\rm otherwise.}\end{array}\right.
$$

Notice that $\bar{G}_{m,k}(\B{Z}_p)\not=\emptyset$  if and only if
$m\le \dfrac{p-1}{2}$, $m$ is even, and $k=\dfrac{m}{2}$. From now
on, we only consider $m\le \dfrac{p-1}{2}$. Let $\sigma\in
S_{\frac{m}{2}}$ and let $\alpha\in\Aut(\B{Z}_p)$. Then
$$|\tilde{F}_{(\alpha,\sigma)}|=\left\{\begin{array}{ll}
\displaystyle \prod_{t=0}^{\frac{m}{2l}-1}
\left(\dfrac{p-1}{2}-rt\right) & \mbox{if
$j_l(\sigma)=\dfrac{m}{2l}$, $|\!<\alpha>\!|=l$,
$l|(\dfrac{p-1}{2})$,
and $l$ is odd,}\\[3ex]
\displaystyle
\prod_{t=0}^{\frac{m}{2l}-1}\left(\dfrac{p-1}{2}-rt\right) &
\mbox{if $j_l(\sigma)=\dfrac{m}{2l}$, $|\!<\alpha>\!|=2l$, and
$l|(\dfrac{p-1}{2})$.}
\end{array}\right.
$$

Notice that $\B{Z}_p$ has no nontrivial subgroup, 
and that the order of each element in $\Aut(\B{Z}_p)$ is a divisor
of $p-1$ and $|\{\alpha\in \Aut(\B{Z}_p)\,:\,
|<\alpha>|=k\}|=\phi(k)$ for each $k|(p-1)$. By summarizing these
together with Theorem~\ref{countcal2}, we find the following
theorem.

\begin{thm}
Let $p$ be a prime number and let $m\le\dfrac{p-1}{2}$. Then
\begin{align*}(p-1)\,\C{E}^w(\B{Z}_p,m)
&=\sum_{k|(\frac{p-1}{2},\frac{m}{2})~and~k~is~odd}\phi(k)\,
\frac{\displaystyle \prod_{t=0}^{\frac{m}{2k}-1}
\left(\frac{p-1}{2}-kt\right)}{\displaystyle
k^{\frac{m}{2k}}\left(\frac{m}{2k}\right)!}\\ &+
\sum_{k|(\frac{p-1}{2},\frac{m}{2})}\phi(2k)\, \frac{\displaystyle
\prod_{t=0}^{\frac{m}{2k}-1}
\left(\frac{p-1}{2}-kt\right)}{\displaystyle
k^{\frac{m}{2k}}\left(\frac{m}{2k}\right)!}.
\end{align*}
In particular, if $\dfrac{m}{2}$ or $\dfrac{p-1}{2}$ is odd, then
$$(p-1)\,\C{E}^w(\B{Z}_p,m) =2\sum_{k|(\frac{p-1}{2},\frac{m}{2})}\phi(k)\,
\frac{\displaystyle \prod_{t=0}^{\frac{m}{2k}-1}
\left(\frac{p-1}{2}-kt\right)}{\displaystyle k^{\frac{m}{2k}}\left(\frac{m}{2k}\right)!}.
 $$
Moreover, if $(\dfrac{p-1}{2},\dfrac{m}{2})=1$, then
$$\C{E}^w(\B{Z}_p,m) =\left(\begin{matrix}\frac12(p-3)\\ \frac{m}{2}\end{matrix}\right).$$ \label{isocayleyzp}
\end{thm}


\begin{thebibliography}{9}

\bibitem{BMPRS} L. Brankovi$\acute{c}$, M. Miller, J. Plesnik, J. Ryan and
J. $\check{S}$ir$\acute{a}\check{n}$, \textit{A note on constructing
large Cayley graphs of given degree and diameter by voltage
assignments}, Electronic Journal of Combinatorics 5 (1998), \#R9.

\bibitem{Cioaba} S. Cioab$\breve{a}$, \textit{Closed walks and eigenvalues of abelian
Cayley graphs}, To appear in C. R. Acad. Sci. Paris, Ser. I.

\bibitem{DS:domonation} I. Dejter and O. Serr, \textit{Efficient dominating sets in Cayley
graphs}, Discrete Applied Mathematics, 129(2) (2003), 319--328.

\bibitem{DSS} C. Droms, B. Servatius and H. Servatius, \textit{Connectivity and
planarity of Cayley Graphs}, Beitr¡§age zur Algebra und Geometrie
Contributions to Algebra and Geometry Volume 39(2) (1998), 269--282.

\bibitem{ET70} B. Elspas and J. Turner, \textit{Graphs with circulant adjacency matrics},
Journal of Combinatorial Theory  9 (1990), 297--307.

\bibitem{FKKL} R. Feng, J. Y. Kim, J. H. Kwak and J. Lee, \textit{Isomorphism classes of
concrete graph coverings}, SIAM J. Discrete Math. 11 (1998),
265--272.

\bibitem{FMT} J. Friedmana, R. Murtyc and J-P. Tillichd, \textit{Spectral estimates for abelian Cayley
graphs}, Journal of Combinatorial Theory, Series B 96 (2006),
111--121.

\bibitem{GTP1} J. L. Gross  and T. W. Tucker, Topological graph theory, Wiley, New York, 1987.

\bibitem{kapo} I. Kapovich, \textit{The geometry of relative Cayley graphs for subgroups of hyperbolic
groups}, preprint, arXiv:math.GR/0201045.

\bibitem{Katznelson} Y. Katznelson, \textit{Chromatic numbers of Cayley graphs on
$\mathbb{Z}$ and recurrence}, Combinatorica 21(2) (2001), 211--219.

\bibitem{KL1} J. H. Kwak and J. Lee, \textit{Isomorphism classes of bipartite cycle
permutation graphs}, ARS Combin. 50 (1998), 139--148.

\bibitem{Li} C. H. Li, \textit{Finite edge-transitive Cayley graphs and rotary
Cayley maps}, Trans. AMS. 358(10) (2006), 4605--4635.

\bibitem{LL} C. H. Li and Z. P. Lu, \textit{Tetravalent edge-transitive Cayley graphs with
odd number of vertices}, Journal of Combinatorial Theory Series B
96(1) (2006), 164--181.

\bibitem{PR} I. Pak and R.
Radoi$\check{c}$i$\acute{c}$, \textit{Hamiltonian paths in Cayley
graphs}, preprint.

\bibitem{Rosenhouse} J. Rosenhouse, \textit{Isoperimetric numbers of Cayley graphs arising from
generalized dihedral groups}, Journal of Combinatorial Mathematics
and Combinatorial Computing 42 (2002), 127--138.

\end{thebibliography}
\end{document}